\documentclass[12pt,twoside,reqno]{article}
\usepackage{amsfonts}
\usepackage{graphicx}
\begin{document}
\begin{center}

{\Large{\textbf{A Generalization of the 2D-DSPM for Solving Linear
System of Equations}}}
\\
\vspace{0.5cm}

{\bf Davod Khojasteh Salkuyeh} \vspace*{0.1cm}

{Department of Mathematics, University of Mohaghegh Ardabili,\\
P. O. Box. 56199-11367, Ardabil, Iran\\
E-mail: khojaste@uma.ac.ir
\\}

\end{center}

\bigskip

\medskip

\begin{center}
 {\bf Abstract}
\end{center}
In  [N. Ujevi$\rm\acute{c}$, New iterative method for solving linear
systems, Appl. Math. Comput. 179 (2006) 725–730], a new iterative
method for solving linear system of equations was presented which
can be considered as a modification of the Gauss-Seidel method. Then
in [Y.-F. Jing and T.-Z. Huang, On a new iterative method for
solving linear systems and comparison results, J. Comput. Appl.
Math., In press]  a different approach, say 2D-DSPM, and more
effective one was introduced. In this paper,  we improve this method
and give a generalization of it. Convergence properties of this kind
of generalization are also discussed. We finally give some numerical
experiments to show the efficiency of the method and compare with
2D-DSPM.


\bigskip

\bigskip

\bigskip

\noindent {\bf\large 1. Introduction}

\bigskip

Consider the linear system of equations
\begin{equation}\label{e:eq01}
    Ax=b,
\end{equation}
where $A\in\mathbb{R}^{n \times n}$ is a symmetric positive definite
(SPD) matrix and $x,b\in \mathbb{R}^{n}$. The Gauss-Seidel method is
an stationary iterative method for solving linear system of equation
and is convergent for SPD matrices.  This method is frequently used
in science and engineering, both for solving linear system of
equations and preconditioning \cite{Benzi,Saadbook}. It can be
easily seen that the Gauss-Seidel method is an special case of a
projection method \cite{Barret,Saadbook}. Let $\mathcal{K}$ and
$\mathcal{L}$ be two $m$-dimensional subspaces of $\mathbb{R}^{n}$.
Let also $x_0$ be an initial guess of the solution of
(\ref{e:eq01}).  A projection method onto $\mathcal{K}$ and
orthogonal to $\mathcal{L}$ is a process which finds an approximate
solution $x\in \mathbb{R}^n$ to (\ref{e:eq01}) by imposing the
conditions that $x\in x_0+\mathcal{K}$ and the new residual vector
be orthogonal to $\mathcal{L}$ (Petrov-Galerkin condition), i.e.
\begin{equation}\label{e:eq02}
    {\rm Find}~~ x \in x_0+\mathcal{K}, \hspace{0.5cm} \textrm{such}~\textrm{that}~\hspace{0.5cm} b-Ax \bot
    \mathcal{L}.
\end{equation}

\noindent It is well known that an iteration of the elementary
Gauss-Seidel method can be viewed as a set of projection methods
with $\mathcal{L}=\mathcal{K}=\{e_i\},~i=1,2,\ldots,n,$ where $e_i$
is the $i$th column of the identity matrix. In fact, a single
correction is made at each step of these projection steps cycled for $i=1,\ldots,n$.\\
\\
\indent In \cite{Ujevic}, Ujevi$\rm\acute{c}$ proposed a
modification of the Gauss-Seidel method which may be named as a
``one-dimensional double successive projection method" and referred
to as 1D-DSPM. In an iteration of 1D-DSPM, a set of double
successive projection methods with two pairs of one-dimensional
subspaces are used. In fact, in an iteration of 1D-DSPM two pairs of
subspaces $(\mathcal{K}_1,\mathcal{L}_1)$ and
$(\mathcal{K}_2,\mathcal{L}_2)$ of one dimension are chosen while it
makes double correction at each step of the process cycled for
$i=1,2,\ldots,n$. In \cite{Jing}, Jing and Huang proposed the
``two-dimensional double successive projection method" and referred
to as 2D-DSPM. In an iteration of 2D-DSPM, a set of projection
methods with a  pairs of two-dimensional subspaces $\mathcal{K}$ and
$\mathcal{L}$ is used and a double correction at each step of the
projection steps is
made.\\
\\
\indent In this paper, a generalization of 2D-DSPM say $m$D-SPM
which is referred to as ``$m$-dimensional successive projection
method" is proposed and its convergence properties are studied.
For $m=2$, the $m$D-SPM  results in 2D-DSPM.\\
\\
\indent Throughout this paper we use some notations as follows. By
$\langle.,.\rangle$ we denote the standard inner product in
$\mathbb{R}^n$. In fact, for two vectors $x$ and $y$ in
$\mathbb{R}^n$, $~\langle x,y\rangle=y^Tx$. For any SPD matrix $M\in
\mathbb{R}^{n \times n}$, the $M$-inner product is defined as
$\langle x,y\rangle_M=\langle Mx,y\rangle$ and its corresponding
norm is  $ \|x\|_M=\langle x,x\rangle_M^{1/2} $. \\
\\
This paper is organized as follows. In section 2, a brief
description of 1D-DSPM and 2D-DSPM are given. In section 3, the
$m$D-SPM is presented and its convergence properties are studied. In
section 4 the new algorithm and its practical implementations are
given. Section 5 is devoted to some numerical experiments to show
the efficiency of the method and comparing with 1D-DSPM and 2D-DSPM.
Some concluding remarks are given in 6.

\bigskip

\bigskip

\medskip

\noindent {\bf\large 2. A brief description of 1D-DSPM and 2D-DSPM}

\bigskip

\medskip

We review 1D-DSPM and 2D-DSPM in the literature of the projection
methods. As we mentioned in the previous section in each iteration
of 1D-DSPM a set of double successive projection method is used. Let
$x_k$ be the current approximate solution. Then the double
successive projection method is applied as following. The first step
is to choose two pairs of the subspaces
$\mathcal{K}_1=\mathcal{L}_1=\{v_1\}$,~
$\mathcal{K}_2=\mathcal{L}_2=\{v_2\}$  and  the next approximate
solution $x_{k+1}$ is computed  as follows
\begin{equation}\label{e:eq03}
    {\rm Find}~~ \widetilde{x}_{k+1}\in x_k+\mathcal{K}_1, \hspace{0.5cm}
    \textrm{such}~\textrm{that}~\hspace{0.5cm} b-A\widetilde{x}_{k+1} \bot
    \mathcal{L}_1,
\end{equation}
\begin{equation}\label{e:eq04}
    {\rm Find}~~ x_{k+1}\in \widetilde{x}_{k+1}+\mathcal{K}_2, \hspace{0.5cm}
    \textrm{such}~\textrm{that}~\hspace{0.5cm} b-Ax_{k+1} \bot
    \mathcal{L}_2.
\end{equation}
This framework results in \cite{Jing,Ujevic}
\[
x_{k+1}=x_k+\alpha_1 v_1+\beta_2 v_2
\]
where
\begin{equation}\label{e:eq05}
\alpha_1=-p_1/a,\qquad \beta_2=(cp_1-ap_2)/ad,
\end{equation}
 in which

\begin{equation}\label{e:eq06}
a=\langle v_1,v_1\rangle_A,\qquad c=\langle v_1,v_2\rangle_A,\qquad
d=\langle v_2,v_2\rangle_A.
\end{equation}
In \cite{Ujevic}, it has been proven that this method is convergent
to the exact solution $x_*$ of
(\ref{e:eq01}) for any initial guess.\\
\\
\indent In the 2D-DSPM, two two-dimensional subspaces
$\mathcal{K}=\mathcal{L}={\rm span}\{v_1,v_2\}$ are chosen and a
projection process onto $\mathcal{K}$ and orthogonal to
$\mathcal{L}$ is used instead of double successive projection method
used in 1D-DSPM. In other words, two subspaces
$\mathcal{K}=\mathcal{L}={\rm span}\{v_1,v_2\}$ are chosen and a
projection method is defined as following.

\begin{equation}\label{e:eq07}
    {\rm Find}~~ x_{k+1}\in x_k+\mathcal{K}, \hspace{0.5cm}
    \textrm{such}~\textrm{that}~\hspace{0.5cm} b-Ax_{k+1} \bot
    \mathcal{L}.
\end{equation}
In \cite{Jing}, it has been shown that this projection process gives
\[
    \alpha=\frac{cp_2-dp_1}{ad-c^2},\qquad     \beta=\frac{cp_1-ap_2}{ad-c^2},
\]
where $p_1,~p_2,~a,~b,~$ and $c$ were defined in Eqs. (\ref{e:eq05})
and (\ref{e:eq06}). It has been proven in \cite{Jing} that the
2D-DSPM is also convergent. Theoretical analysis and numerical
experiments presented in \cite{Jing} show that 2D-DSPM is more
effective than the 1D-DSPM.\\
\\
\indent A main problem with  1D-DSPM and 2D-DSPM is to choose the
optimal vectors $v_1$ and $v_2$. In this paper, we first propose a
generalization of 2D-DSPM and then give a strategy to choose vectors
$v_1$ and $v_2$ in a special case.

\bigskip

\bigskip

\medskip

\noindent {\bf\large 3. $m$-dimensional successive projection
method}

\bigskip

\medskip

Let $\{v_1,v_2,\ldots,v_m\}$ be a set of $m$ independent vectors in
$\mathbb{R}^n$. For later use, let also $V_m=[v_1,v_2,\ldots,v_m]$.
Now we define the $m$D-SPM as follows. In an iteration of the
$m$D-SPM we use a set of projection process onto $\mathcal{K}={\rm
span}\{v_1,v_2,\ldots,v_m\}$ and orthogonal to
$\mathcal{L}=\mathcal{K}$. In other words, two $m$-dimensional
subspaces $\mathcal{K}$ and $\mathcal{L}$ are used in the projection
step instead of two two-dimensional subspaces used in 2D-DSPM. In
this case Eq. (\ref{e:eq02}) turns the form
\begin{equation}\label{e:eq08}
   {\rm Find}~~y_m\in \mathbb{R}^m \quad \textrm{such}~\textrm{that}
\quad x_{k+1}=x_k+V_my_m, \quad {\rm and} \quad V_m^T(
b-Ax_{k+1})=0.
\end{equation}
We have
\begin{eqnarray*}
  r_{k+1} &=& b-Ax_{k+1} \\
          &=& b-A(x_k+V_my_m)\\
          &=& r_k-AV_my_m,
\end{eqnarray*}
where $r_k=b-Ax_k$. Hence from Eq. (\ref{e:eq08}) we deduce
\[
0=V_m^Tr_{k+1}=V_m^T(r_k-AV_my)=V_m^Tr_k-V_m^TAV_my_m \Rightarrow
V_m^TAV_my_m=V_m^Tr_k.
\]
The matrix $V_m^TAV_m$ is an SPD matrix, since $A$ is SPD. Therefore
\begin{equation}\label{e:eq09}
    y_m=(V_m^TAV_m)^{-1}V_m^Tr_k.
\end{equation}
Hence, from (\ref{e:eq08}) we conclude that
\begin{equation}\label{e:eq10}
    x_{k+1}=x_k+V_m(V_m^TAV_m)^{-1}V_m^Tr_k.
\end{equation}

\bigskip

\noindent {\bf Theorem 1.} {\it Let $A$ be an SPD matrix and assume
that $x_k$ is an approximate solution of (\ref{e:eq01}). Then
\begin{equation}\label{e:eq11}
\| d_k \|_A \geq \|d_{k+1}\|_A,\\
\end{equation}
where $d_k=x_*-x_k$ and $d_{k+1}=x_*-x_{k+1}$ in which $x_{k+1}$ is
the approximate solution computed by Eq. (\ref{e:eq10}).}

\bigskip

\noindent {\bf Proof.} It can be easily verified that
$d_{k+1}=d_k-V_my_m$ and $Ad_k=r_k$ where $y_m$ is defined by
(\ref{e:eq07}). Then
\begin{eqnarray*}
\langle Ad_{k+1},d_{k+1} \rangle& = & \langle Ad_k-AV_my_m,d_k-V_my_m \rangle \\
& =  & \langle Ad_k,d_k \rangle-2 \langle V_my_m,r_k \rangle
       + \langle  AV_my_m,V_my_m \rangle \qquad {\rm by}~ Ad_k=r_k  \\
& =  & \langle Ad_k,d_k \rangle - \langle y_m,V_m^Tr_k \rangle
\hspace{4.4cm} {\rm by}~(\ref{e:eq09}).
\end{eqnarray*}
Therefore
\begin{eqnarray*}
  \| d_k \|_A^2 - \|d_{k+1}\|_A^2 &=& \langle Ad_k,d_k \rangle-\langle Ad_{k+1},d_{k+1} \rangle \\
                                  &=& \langle y_m,V_m^Tr_k \rangle\\
                                  &=& (V_m^Tr_k)^T (V_m^TAV_m)^{-1} V_m^Tr_k,  \qquad {\rm
                                     by}~(\ref{e:eq09})\\
                                  &=:& S(r_k).
\end{eqnarray*}
Since $(V_m^TAV_m)^{-1}$ is SPD then we have
\[
\| d_k \|_A^2 - \|d_{k+1}\|_A^2 \geq 0,
\]
and the desired result is obtained. \qquad\qquad $\Box$ \\
\\
\indent This theorem shows that if $V_m^Tr_k=0$ then $S(r_k)=0$ and
we don't have any reduction in the square of the $A$-norm of  error.
But, if $V_m^Tr_k \neq 0$ then the square of the $A$-norm of error
is reduced by $S(r_k)>0)$.
\\
\\
\noindent {\bf Theorem 2.} {\it Assume that $A$ is an SPD matrix and
$\mathcal{L}=\mathcal{K}$. Then a vector $x_{k+1}$  is the result of
projection method onto $\mathcal{K}$ orthogonal to $\mathcal{L}$
with the starting vector $x_k$ iff it minimizes the $A$-norm of the
error over $x_k+\mathcal{K}$.}

\bigskip

\noindent {\bf Proof.} See \cite{Saadbook}, page 126. \qquad\qquad $\Box$ \\

This theorem shows that if $\mathcal{K}_1=\mathcal{L}_1 \subset
\mathcal{K}_2=\mathcal{L}_2$, then the reduction of $A$-norm of the
errors obtained by the subspaces $\mathcal{K}_2=\mathcal{L}_2$ is
more than or equal to that of the subspaces
$\mathcal{K}_1=\mathcal{L}_1$. Hence by increasing the value of $m$,
the convergence rate may increase.\\
\\
\indent In the continue we consider the special case that the
vectors $v_i$ are the column vectors of the identity matrix. The
next theorem not only proves the convergence of the method in this
special case but also gives an idea to choose the optimal vectors
$v_i$.\\
\\
\noindent {\bf Theorem 3.} {\it Let $\{i_1,i_2,\ldots,i_m \}$ be the
set of indices of $m$ components of largest absolute values in $r_k$
such that $i_1<i_2<\ldots < i_m$. If $v_j=e_{i_j},~j=1,\ldots,m$
then
\begin{equation}\label{e:eq12}
\| d \|_A^2 - \|d_{new}\|_A^2 \geq \frac{m}{n \lambda_{\max}(A)
}\|r_k\|_2^2.
\end{equation}}
\bigskip
\noindent {\bf Proof.} Let $E_m=[e_{i_1},e_{i_2},\ldots,e_{i_m}]$.
By Theorem 1, we have
\begin{eqnarray*}
  \| d_k \|_A^2 - \|d_{k+1}\|_A^2 &=& S(r_k) \\
          &=& (E_m^Tr_k)^T (E_m^TAE_m)^{-1}E_m^Tr_k.
\end{eqnarray*}
Then by using Theorem 1.19 in \cite{Saadbook} we conclude

\begin{equation}\label{e:eq13}
S(r_k) \geq \lambda_{\min}((E_m^T A E_m)^{-1}) \|E_m^Tr_k\|_2^2 \geq
\frac{1}{\lambda_{\max}(E_m^T A E_m)} \|E_m^T r_k \|_2^2,
\end{equation}
where for a square matrix $Z$, $\lambda_{\min}(Z)$ and
$\lambda_{\max}(Z)$ stand for the smallest and largest eigenvalues
of $Z$. It can be easily verified that \cite{DKS}

\begin{equation}\label{e:eq14}
\lambda_{\max} (E_m^T A E_m)\leq \lambda_{\max}(A), \quad
\frac{(E_m^T r_k,E_m^Tr_k)}{(r_k,r_k)} \geq \frac{m}{n}.
\end{equation}
Hence
\begin{equation}\label{e:eq15}
    S(r_k) \geq \frac{m}{n \lambda_{\max}(A)} \|r_k\|_2^2,
\end{equation}
and the desired result is obtained. \qquad\qquad $\Box$\\
\\
\indent Eq. (12) shows the convergence of the method. Eq.
(\ref{e:eq13}) together with the first relation of the equation
(\ref{e:eq14}) give
\[
S(r_k) \geq  \frac{1}{\lambda_{\max}(A)} \|E_m^T r_k \|_2^2.
\]
This equation gives an idea to choose indices $i_j,~j=1,\ldots,m$.
In fact, it shows that if these indices are the $m$ components of
the largest absolute values in $r_k$, then the lower bound of
$S(r_k)$ depends on $\|E_m^T r_k \|_2^2$, and will be as large as
possible.\\
\\
\indent In \cite{DKS}, another theorem for the convergence of the
method obtained by $v_i=e_{i_j}$ was presented and an algorithm
based upon this theorem was constructed for computing a sparse
approximate inverse factor of an SPD matrix and was used as a
preconditioner for SPD linear system.
\bigskip

\bigskip

\medskip

\noindent {\bf\large 4. Algorithm and its practical implementations}

\bigskip

\medskip

\indent Hence, according to the results obtained in the previous
section we can summarized the $m$D-DSM in the special case that
$v_j=e_{i_j}$ as following.
\bigskip

\medskip

\underline{{\bf Algorithm 1:} $m$D-DSM}\\
\vspace{-0.5cm}

\begin{itemize}
\item[1.] Choose an initial guess $x_0\in\mathbb{R}^n$ to (\ref{e:eq01}) and $r=b-Ax_0$.\vspace{-0.3cm}

\item[2.] Until convergence, Do \vspace{-0.3cm}

\item[3.] \qquad $x:=x_0$ \vspace{-0.3cm}

\item[4.] \qquad For $i=1,\ldots,n$,  Do \vspace{-0.3cm}

\item[5.] \qquad \qquad Select the indices $i_1,i_2,\ldots,i_m$ of $r$ as defined in Theorem
3 \vspace{-0.3cm}

\item[6.] \qquad \qquad $E_m:=[e_{i_1},e_{i_2},\ldots,e_{i_m}]$ \vspace{-0.3cm}

\item[7.] \qquad \qquad Solve $(E_m^TAE_m)y_m=E_m^Tr$ for $y_m$ \vspace{-0.3cm}

\item[8.] \qquad \qquad $x:=x+E_my_m$ \vspace{-0.3cm}

\item[9.]\qquad \qquad $r:=r-AE_my_m$ \vspace{-0.3cm}

\item[10.] \qquad EndDo \vspace{-0.3cm}

\item[11.] \qquad $x_0:=x$ and if $x_0$ has converged then Stop  \vspace{-0.3cm}

\item[12.] EndDo
\end{itemize}

\bigskip

In practice, we see that the matrix $E_m^TAE_m$ is a principal
submatrix of $A$ with column and row indices  in
$\{i_1,i_2,\ldots,i_m\}$. Hence, we do not need any computation for
computing the matrix $E_m^TAE_m$ in step 7. For solving the linear
system in step 7 of this algorithm one can use the Cholesky
factorization of the coefficient matrix. Step 8 of the algorithm may
be written as

\begin{itemize}
  \item For $j=1,2,\ldots,m$\vspace{-0.2cm}
  \item \qquad $x_{i_j}:=x_{i_j}+(y_m)_j$\vspace{-0.2cm}
  \item EndDo
\end{itemize}
Hence only $m$ components of the vector $x$ are modified. Step 9 of
the algorithm can be written as
\[
r=r-\sum_{j=1}^{m}(y_m)_j ~ a_{i_j}
\]
where $a_{i_j}$ is the $i_j$ column of the matrix $A$.\\
\\
\indent It can be seen that Algorithm 2 in \cite{Jing} is an special
case of this algorithm. In fact, if $m=2$ and the indices $i_1$ and
$i_2$ are chosen as $i_1=i$ and $i_2=i-ij_{gap}$ ($i_2=i-ij_{gap}+n$
if $i\leq ij_{gap}$), where $ij_{gap}$ is a positive
integer parameter less than $n$, then Algorithm 2 in
\cite{Jing} is obtained.\\
\\
\indent As we see, the first advantage of our algorithm over
Algorithm 2 in \cite{Jing} is that our algorithm chooses the indices
$\{i_1,i_2,\ldots,i_m\}$, automatically. Another advantage is that
our algorithm chooses the indices such that the reduction in the
square of the $A$-norm of the error is more than that of Algorithm 2
in \cite{Jing}. Numerical experiments in the next section also
confirm also this fact. The main advantage of Algorithm 2 in
\cite{Jing} over our algorithm is that only $m$ components of the
current residual are computed whereas in our algorithm the residual
vector should be computed for choosing $m$ indices of largest
components in absolute value.

\bigskip

\bigskip

\medskip

\noindent {\bf 5. Numerical experiments}

\bigskip

\medskip

In this section we give some numerical experiments to compare our
method with Algorithm 2 in \cite{Jing}. Numerical results have been
obtained by some MATLAB codes. We use all of the assumptions such as
initial guess, exact solution, stopping criterion, and the examples
used in \cite{Jing}. Let $b=Ae$, where $e$ is an $n$-vector whose
elements are all equal to unity, i.e., $~e=(1,1,\ldots,1)^T$. We use
$\|x_{k+1}-x_{k}\|_{\infty} < 10^{-6}$ as the stopping criterion. An initial
guess equal to $x_0=(x^1,\ldots,x^n)$, where $x^i=0.001 \times i$,
$i=1,2,\ldots,n$ is chosen. For each of the systems we give the
numerical experiments of Algorithm 2 in \cite{Jing} with
$ij_{gap}=2$ and $500$ and our algorithm with $m=2,3,4$ and $5$.\\
\\
\noindent {\bf Example 1.} Let $A=(a_{ij})$ where
\[
a_{ii}=4n, \quad a_{i,i+1}=a_{i+1,i}=n, \quad a_{ij}=0.5 \quad {\rm
for}~i=1,2,\ldots,n,~ j \neq i,i+1.
\]
We also assume $n=1000$. Numerical experiments in terms of iteration
number were shown in Table 1.

\begin{table}
\centering \caption{Results for the Example
1}\label{T1}\vspace{0.3cm}

\begin{tabular}{l c c c c}

\hline

Algorithm 2 in \cite{Jing} & 2D-DSM   &  3D-DSM   & 4D-DSM &
5D-DSM\\ \hline

6~ ($ij_{gap}=2$)  & 5  & 4 & 3 & 2\\

7~ ($ij_{gap}=500$)& \\ \hline
\end{tabular}
\end{table}

Numerical experiments presented in Table 1 show that the 2D-DSM
method gives better results than the Algorithm 2 in \cite{Jing}.
This table also shows the effect of increase in $m$ on the number of
iterations for convergence.\\
\\
\noindent {\bf Example 2.} Let $A=(a_{ij})$ be the same matrix used
in the previous example except the diagonal entries are changed to
\[
a_{ii}=3n,~ i=1,2,\ldots,n.
\]
Numerical experiments were given in Table 2.

\begin{table}
\centering \caption{Results for the Example
2}\label{T2}\vspace{0.3cm}

\begin{tabular}{l c c c c}

\hline

Algorithm 2 in \cite{Jing} & 2D-DSM   &  3D-DSM   & 4D-DSM &
5D-DSM\\ \hline

8~ ($ij_{gap}=2$)  & 7  & 6 & 4 & 4\\

9~ ($ij_{gap}=500$)& \\ \hline
\end{tabular}
\end{table}
This table also shows the advantages of our method on Algorithm 2
\cite{Jing}.\\
\\
\noindent {\bf Example 3.} Our third set of test matrices used arise
from standard five point finite difference scheme to  discretize
\[
-\triangle u+a(x,y)u_x+b(x,y)u_y + c(x,y ) u = f(x,y),
\quad {\rm in} \quad \Omega=[0,1] \times [0,1], \\
\]
where $a(x,y),~b(x,y),~c(x,y)$ and $d(x,y)$ are given real valued
functions. We consider three following cases:

\noindent$\quad {\rm Case 1}:
a(x,y)=0,~~  b(x,y)=10(x+y),~~ c(x,y)=10(x-y),~~ f(x,y)=0,\\
~\quad {\rm Case 2}: a(x,y)=-10(x+y),~~  b(x,y)=-10(x-y),~~
c(x,y)=1,~~
f(x,y)=0,\\
~\quad{\rm Case 3}: a(x,y)=10e^{xy},~~  b(x,y)=10e^{-xy},~~
c(x,y)=0,~~
f(x,y)=0.$\\
We assume $m=32$. In this case we obtain three SPD matrices of order
$n=32 \times 32$ \cite{Jing} and used them as the coefficient of the
linear systems. Numerical results were given in Table 3.

\begin{table}
\centering \caption{Results for the Example
3}\label{T3}\vspace{0.3cm}

\begin{tabular}{l l c c c c}

\hline

 Cases & Algorithm 2 in \cite{Jing} & 2D-DSM   &  3D-DSM   & 4D-DSM &
5D-DSM\\ \hline

Case 1 & 391 ~ ($ij_{gap}=2$)  & 226  & 153 & 116 & 94\\

       & 323 ~ ($ij_{gap}=500$)& \\
\\
Case 2 & 312 ~ ($ij_{gap}=2$)  & 192  & 131 & 100 & 80\\

       & 256 ~ ($ij_{gap}=500$)& \\
\\
Case 3 & 302 ~ ($ij_{gap}=2$)  & 218  & 151 & 115 & 93\\

       & 250~ ($ij_{gap}=500$)& \\
\hline
\end{tabular}
\end{table}
This table also confirm that our method is more effective that the
Algorithm 2 \cite{Jing}.

\bigskip

\bigskip

\medskip

\noindent {\bf 5. Conclusion}

\bigskip

In this paper a generalization of the 2D-DSPM \cite{Jing} which
itself is a generalization of 1D-DSPM \cite{Ujevic} is presented.
1D-DSPM and 2D-DPM need to prescribed some subspaces of
$\mathbb{R}^n$ for the projection steps. But our method in the
spacial case chooses this subspaces automatically. Theoretical
analysis and numerical experiments presented in this paper showed
that our method is more effective that 2D-DSPM.

\bigskip

\bigskip

\medskip

\noindent {\bf 6. Acknowledgments}

\bigskip

The author would like to thank Yan-Fei Jing for providing the
matrices of Example 3.

\end{document}